\documentclass{amsart}
\usepackage{mathrsfs}
\usepackage{amsfonts,amsmath,amscd,amssymb,amsbsy,amsthm,amstext,amsopn,amsxtra,mathtools}
\usepackage{color,mathrsfs,verbatim,subfigure}
\usepackage{extarrows}
\usepackage[all]{xy}
\usepackage{caption}
\usepackage{enumitem}
\usepackage{pdfpages}
\usepackage[toc,page]{appendix}
\usepackage[T1]{fontenc}
\usepackage{amsmath}
\usepackage{stmaryrd,amssymb}
\usepackage[utf8x]{inputenc}
\usepackage{t1enc}

\newtheorem{thm}{Theorem}[section]

\newtheorem{cor}[thm]{Corollary}
\newtheorem{prop}[thm]{Proposition}

\newtheorem{defn}{Definition}

\newtheorem{rmk}{Remark}

\newtheorem{example}{Example}

\renewcommand{\Re}{\operatorname{Re}}
\renewcommand{\Im}{\operatorname{Im}}

\newcommand{\CC}{\mathbb{C}}

\newcommand{\RR}{\mathbb{R}}
\newcommand{\ZZ}{\mathbb{Z}}

\newcommand{\cV}{\mathcal{V}}
\newcommand{\cX}{\mathcal{X}}

\newcommand{\fM}{\mathfrak{M}}

\newcommand{\Def}{\operatorname{Def}}

\newcommand{\divv}{\operatorname{div}}

\newcommand{\Kt}{\operatorname{K}}
\newcommand{\Kum}{\operatorname{Kum}}
\newcommand{\Mon}{\operatorname{Mon}}

\newcommand{\Or}{\operatorname{O}}

\newcommand{\Pic}{\operatorname{Pic}}

\newcommand{\rk}{\operatorname{rk}}

\title[Connected components of manifolds of Kummer type]{Connected components of moduli spaces of irreducible holomorphic symplectic manifolds of Kummer type}
\author{Claudio Onorati}
\address{Department of Mathematical Sciences\\University of Bath\\Bath BA2 7AY\\England} 
\curraddr{Dipartimento di Matematica, Universit\`a degli studi di Roma Tor Vergata, via della Ricerca Scientifica 1, 00133 Roma, Italy}
\email{onorati@mat.uniroma2.it}
\begin{document}
\maketitle

\begin{abstract} 
In this paper we determine the number of connected components of moduli spaces
of both marked and polarised irreducible holomorphic symplectic manifolds deformation
equivalent to generalised Kummer varieties. 
\end{abstract}

\section{Introduction} 

An irreducible holomorphic symplectic manifold $X$ is a simply connected compact K\"ahler manifold with an everywhere nondegenerate holomorphic $2$-form that is unique up to constant.
They are seen as higher dimensional analogous to $K3$ surfaces, that are the lowest dimensional 
example, and with which they share a lot of properties. 
The existence of the nondegenerate form implies that the dimension 
is even. For any $n\geq2$, Beauville constructed two families of irreducible holomorphic
symplectic manifolds of dimension $2n$, which are not deformation equivalent each other 
(\cite{Beauville:c1=0}). These families are constructed starting from Hilbert schemes of
points on $K3$ and abelian surfaces, and are respectively called Hilbert powers of $K3$
surfaces and generalised Kummer varieties. 
Together with two sporadic families constructed by O'Grady
in dimension $6$ and $10$ (see \cite{O'Grady:6dim} and \cite{O'Grady:10dimPublished}), 
and not deformation equivalent to the previous ones, these are 
all the known deformation types of irreducible holomorphic symplectic manifolds. 
The question whether this is a complete list or not is still open.
In this paper we investigate some topological properties of moduli spaces of one of these types, 
namely the generalised Kummer deformation type, and we simply refer to these manifolds as \emph{manifolds of Kummer type}.

More precisely we focus on moduli problems. Fixed a deformation type, one can costruct moduli 
spaces of marked and polarised irreducible holomorphic symplectic manifolds. These moduli
spaces behave very much like the corresponding moduli spaces of $K3$ surfaces. For example,
we have local and global Torelli theorems that create a bridge between geometry and
combinatorics (see for example \cite{GrossHuybrechtsJoyce:CalabiYau} and 
\cite{Markman:Survey}). We address the topological question: how many connected
components do these moduli spaces have?  

For $K3$ surfaces it is known that the marked moduli space has
two connected components, corresponding to the two pairs $(S,\eta)$ and $(S,-\eta)$,
where $S$ is a $K3$ surface and $\eta$ is a marking. In other
words, this reflects the fact that the monodromy group $\Mon^2(K3)$ has index $2$
in the isometry group $\Or(H^2(K3,\ZZ))$, and $-id$ is not a monodromy operator.
On the other hand, fixing a polarisation does
not add connected components, i.e.\ the moduli space of polarised $K3$ surfaces is
connected (this is mostly a consequence of the fact that the $K3$ lattice $H^2(K3,\ZZ)$
is unimodular). It is quite natural to ask what is the situation in higher dimensions.

Let us also point out that determining the number of connected components also gives
information about the geometry itself. 
If $X$ and $X'$ are two non-birational irreducible holomorphic symplectic manifolds, that are
Hodge isometric (i.e.\ they have the same period), then there exist markings $\eta$ and
$\eta'$ such that $(X,\eta)$ and $(X',\eta')$ belong to different connected components.
This is a consequence of the Hodge-theoretic Torelli theorem (\cite[Theorem~1.3]{Markman:Survey}).
For example, Namikawa noticed that the generalised Kummer varieties constructed on an abelian 
surface and its dual
are not birational, at least when $A$ is generic, but they are Hodge isometric (\cite{Namikawa:Counterexample}).
Markman and Mehrotra also constructed a Hodge isometry between this two varieties in
\cite{MarkmanMehrotra:HilbertSchemeDense}. As a consequence one deduces
that the moduli space of marked manifolds of Kummer type 
has always at least four connected components, and this bound is reached in the special
situation when $n+1$ is the power of a prime (recall that $2n$ is the dimension of the manifold). 
Even though similar phenomena appear for Hilbert powers of $K3$ surfaces as well in certain 
dimensions (cf.\ \cite[Proposition~4.10]{Markman:IntegralConstraints}), when $n-1$ is the power of a prime the moduli space of marked manifolds of $K3^{[n]}$ type has only two connected components. This difference between Hilbert powers of $K3$ surfaces and generalised Kummer varieties is due to the 
aforementioned Namikawa phenomenon, which is very geometrical in nature.

As we recall in Section~\ref{section:pre}, the set of connected components is strictly
related to the shape of the monodromy group. Monodromy groups are important gadgets
attached to any irreducible holomorphic symplectic manifolds, that can be naively thought of as
groups of geometric isometries of the $H^2$ lattice 
(see Definition~\ref{defn:parallel transport operators}). Their computation is paramount to
study the geometry of irreducible holomorphic symplectic manifolds, and they have now be computed in all the known deformation classes (see \cite{Markman:Monodromy1} for the $K3^{[n]}$ type, \cite{Mongardi:Monodromy} and \cite{Markman:MonodromyKum} for the Kummer type, \cite{MongardiRapagnetta:MonodromyOG6} for the OG6 type, and \cite{Onorati:Monodromy} for the OG10 type). %So far, monodromy groups are known only for Hilbert powers of $K3$ surfaces (\cite{Markman:Monodromy1}) and generalised Kummer varieties (\cite{Mongardi:Monodromy}, \cite{Markman:MonodromyKum}). 

Apostolov in \cite{Apostolov:ConnectedComponents} computes the number of connected
components of moduli spaces of both marked and polarised irreducible holomorphic symplectic
manifolds deformation equivalent to Hilbert powers of $K3$ surfaces. We extend his results
to manifolds of Kummer type. 
Since the Beauville-Bogomolov-Fujiki lattice of Hilbert powers of $K3$ surfaces and generalised
Kummer varieties are both of the form $L\oplus\langle l\rangle$, where $L$ is unimodular
and $l^2<0$, our result is, as expected, very similar to Apostolov's result. 

The main result of the paper is the following.

\begin{thm}
Let $n\geq2$.
\begin{enumerate}
\item ({\bf Corollary~\ref{cor:marked}}) The moduli space of marked irreducible holomorphic symplectic 
manifolds of Kummer type has $2^{\rho(n+1)+1}$
connected components, where $\rho(k)$ is the number of distinct primes in the factorisation 
of $k$.

\item ({\bf Theorem~\ref{thm:connected components polarised}}) The moduli space of polarised 
irreducible holomorphic symplectic manifolds of Kummer type is not connected in general. Moreover, the number of connected components can get 
arbitrarily large as the degree of the polarisation increases 
(see Theorem~\ref{thm:connected components polarised} for the precise statement).
\end{enumerate}
\end{thm}
%\vspace{0.3cm}

The main tool used to prove the Theorem above is a characterisation of (polarised)
parallel transport operators in terms of a distinguished orbit of primitive embeddings of the
Beauville-Bogomolov-Fujiki lattice in the Mukai lattice 
(see Proposition~\ref{prop:characterisation pto} and the Corollary thereafter). This
characterisation solves \cite[Problem~10.3]{Markman:Survey} for manifolds of Kummer type.
\vspace{0.3cm}

{\bf Plan of the paper.} In Section~\ref{section:pre} we recall the main definitions,
notations and background to state and prove the results in the rest of the paper.
In Section~\ref{section:pto} we recall the definition of generalised Kummer varieties
and give a characterisation of (polarised) parallel transport operators 
(cf.\ Proposition~\ref{prop:characterisation pto}). In Section~\ref{section:marked} we 
focus on the moduli space of marked irreducible holomorphic symplectic manifolds of 
Kummer type; in Section~\ref{section:polarised} we focus on the moduli space of polarised 
irreducible holomorphic symplectic manifolds of Kummer type.
\vspace{0.3cm}

{\bf Acknowledgments.} I would like to thank Professor Gregory Sankaran for the great amount of conversations and advice. Moreover I thank Giovanni Mongardi for having read a preliminary version of this work and Benjamin Wieneck for communicating his work to me. Finally I thank the anonymous referee for the many comments on the paper, which have improved it a lot.
\vspace{0.3cm}

{\bf Note added in proof.} This manuscript has lived a long time as a preprint before being submitted for publication. In the meantime progresses have been done in the understanding of the other two known deformation types, namely the so-called OG6 and OG10 deformation types. It has been proved in \cite{MongardiRapagnetta:MonodromyOG6} and \cite{Onorati:Monodromy} that in both these cases the monodromy group is maximal, i.e.\ it coincides with the group of orientation preserving isometries. In particular the number of connected components is easy to compute in these cases: we have two connected components for the moduli space of marked manifolds of type OG6 and OG10 (corresponding to change the sign of the marking), and only one connected component for the moduli space of polarised manifolds with fixed polarisation type.

%%%%%%%%%%%%%%%%%%%%%%%%%%%%%%%%%%%%%%%%%
\section{Preliminaries and notations}\label{section:pre}
\begin{defn}
A compact K\"ahler manifold $X$ is called \emph{irreducible holomorphic symplectic} if it is 
simply connected and $H^0(X,\Omega^2_X)=\CC\sigma_X$, where $\sigma_X$ is nondegenerate
at any point.
\end{defn} 
It follows directly from the definition that $H^2(X,\ZZ)$ is a torsion free $\ZZ$-module;
it turns out to be a lattice thanks to the Beauville-Bogomolov-Fujiki form $q_X$
(see for example \cite{Beauville:c1=0}). This lattice structure is paramount to study the 
geometry of an irreducible holomorphic symplectic manifold $X$; we refer to 
\cite{GrossHuybrechtsJoyce:CalabiYau} and \cite{Markman:Survey} for a detailed account
of results on their geometry.

Let $X_1$ and $X_2$ be two irreducible holomorphic symplectic manifolds that are
deformation equivalent.
\begin{defn}\label{defn:parallel transport operators}
We say that $g\colon H^2(X_1,\ZZ)\to H^2(X_2,\ZZ)$ is a \emph{parallel transport operator} 
if there exists a family $p\colon\cX\to B$, points $b_1,b_2\in B$ and isomorphisms 
$\varphi_i\colon X_i\stackrel{\sim}{\longrightarrow}\cX_{b_i}$ such that the composition 
$(\varphi_2^*)^{-1}\circ g\circ\varphi_1^*$ is the parallel transport inside the local system 
$R^2 p_*\ZZ$ along a path $\gamma$ from $b_1$ to $b_2$. 
Here $R^2p_*\ZZ$ is endowed with the
Gauss-Manin connection.

When $X_1=X_2=X$ and $\gamma$ is a loop we talk about \emph{monodromy operators}. 
Such isometries form a subgroup $\Mon^2(X)$ of the orthogonal group $\Or(H^2(X,\ZZ))$, 
called \emph{monodromy group}.
\end{defn}
\begin{rmk}\label{rmk:orient pub}
If $X$ is an irreducible holomorphic symplectic manifold, the Beauville--Bogomolov--Fujiki lattice $H^2(X,\ZZ)$ has signature $(3,b_2(X)-3)$. As shown in \cite[Lemma~4.1]{Markman:Survey}, the cone $\widetilde{C}_X$ of positive classes $H^2(X,\RR)$ (not to be confused with the positive cone of $X$) has a $1$-dimensional cohomology space, i.e.\ $H^2(\widetilde{C}_X,\ZZ)\cong\ZZ$. A choice of a generator is usually referred to as an \emph{orientation} of $H^2(X,\ZZ)$. Any isometry of $H^2(X,\ZZ)$ induces an automorphism of $H^2(\widetilde{C}_X,\ZZ)$, and an isometry is said to be \emph{orientation preserving} if it is in the kernel of the map $\Or(H^2(X,\ZZ))\to\operatorname{Aut}(H^2(\widetilde{C}_X,\ZZ))\cong\ZZ_2$. By a direct check, or again by \cite[Lemma~4.1]{Markman:Survey}, reflections around vectors of degree $-2$ are orientation preserving; reflections around vectors of degree $2$ are not orientation preserving. We remark that, up to a sign, this definition of orientation preserving isometry coincides with the notion of spinor norm in lattice theory. 

If now $\omega_X$ is a K\"ahler class and $\sigma_X$ is a holomorphic symplectic $2$-form, then the positive $3$-space $\langle\Re(\sigma_X),\Im(\sigma_X),\omega_X\rangle\subset
H^2(X,\RR)$ determines a preferred orientation: in fact, by \cite[Lemma~4.1]{Markman:Survey} again, for any positive $3$-space $W$ in $H^2(X,\RR)$, the space $W\setminus0$ is a retract of $\widetilde{C}_X$. Notice that the orientation so defined does not depend on the choice of the K\"ahler class, nor on the choice of the symplectic form (see \cite[Section~4]{Markman:Survey}). In particular, if $g\colon H^2(X,\ZZ)\to H^2(Y,\ZZ)$ is an isometry, then we say that $g$ is \emph{orientation preserving} if it preserves the preferred orientations of $X$ and $Y$. 

Notice that, by definition, any parallel transport operator is orientation preserving.
If $\Or^+(H^2(X,\ZZ))$ denotes the group of orientation 
preserving isometries, then $\Mon^2(X)\subset\Or^+(H^2(X,\ZZ))$.

%For any irreducible holomorphic symplectic manifold $X$, let us denote by $\omega_X$ the K\"ahler class and by $\sigma_X$ the symplectic form. The positive (real) three-space $\langle\Re(\sigma_X),\Im(\sigma_X),\omega_X\rangle\subset H^2(X,\RR)$ comes then with a preferred orientation. We say that an isometry $H^2(X_1,\ZZ)\to H^2(X_2,\ZZ)$ is \emph{orientation preserving} if it preserves this orientation. By definition, any parallel transport operator is orientation preserving. In particular, if $\Or^+(H^2(X,\ZZ))$ denotes the group of orientation preserving isometries, then $\Mon^2(X)\subset\Or^+(H^2(X,\ZZ))$.
\end{rmk}
Let $\Lambda$ be a lattice. 
A \emph{marking} is an isometry $\eta\colon H^2(X,\ZZ)\to\Lambda$. We denote by
$\fM_\Lambda$ the moduli space parametrising pairs $(X,\eta)$ where $X$ has a fixed
deformation type, and $\eta$ is a marking. 
$\fM_\Lambda$ is a smooth (but not Hausdorff) complex manifold of dimension 
$\rk\Lambda-2$ (complex charts are given by Kuranishi spaces $\Def(X)$).
Let $\Xi_\Lambda$ be the set of connected components of $\fM_\Lambda$.

The group $\Or(\Lambda)$ acts on $\fM_\Lambda$ by changing the marking,
and the induced action on $\Xi_\Lambda$ is transitive. Moreover, the stabiliser of a connected 
component is a group isomorphic to the monodromy group (see \cite{Markman:Survey} for
the detailed statement).
Therefore, the cardinality of $\Xi_\Lambda$ is equal to the index of the monodromy group
$\Mon^2(X)$ in $\Or(H^2(X,\ZZ))$, where $X$ is any irreducible holomorphic symplectic
manifold in $\mathfrak{M}_\Lambda$. By \cite[Theorem~2.6]{Huybrechts:Finiteness}
and \cite[Theorem~3.5.(iv), Theorem~7.2]{Verbitsky:GTT}, this number is finite.  

Now, let $h\in\Lambda$ be a primitive and positive element; we denote by $\bar{h}$ the 
$\Or(\Lambda)$-orbit determined by $h$. The moduli space $\mathfrak{M}^a_{\bar{h}}$
parametrising triples $(X,H,\eta)$, where $\eta$ is a marking and $H$ is a polarisation on $X$
such that $\eta(c_1(H))\in\bar{h}$, is the moduli space of marked polarised irreducible
holomorphic symplectic manifolds of type $\bar{h}$ (cf.\ \cite[Section~7]{Markman:Survey}). 
If we denote by $\mathfrak{M}^{t,a}_{\bar{h}}$ a connected component of $\mathfrak{M}^a_{\bar{h}}$, then by \cite[Lemma~8.1]{Markman:Survey} there is an isomorphism of analytic spaces 
\[ \mathfrak{M}^a_{\bar{h}}/\Or(\Lambda)\cong \mathfrak{M}^{t,a}_{\bar{h}}/\Gamma, \]
where $\Gamma$ is an arithmetic subgroup of $\Or(\Lambda)$ isomorphic to the subgroup of the monodromy group preserving the polarisation.

By \cite[Lemma~8.3]{Markman:Survey} 
(cf.\ also \cite[Theorem~1.5]{GritsenkoHulekSankaran:Moduli}), the quotient 
$\cV^t_{\bar{h}}:=\mathfrak{M}^{t,a}_{\bar{h}}/\Gamma\cong\mathfrak{M}^a_{\bar{h}}/\Or(\Lambda)$ is (analytically) isomorphic to a connected component 
of the moduli space $\cV_{n,d}$ of polarised irreducible holomorphic symplectic manifolds of dimension $2n$ and
with a polarisation of degree $d$. Here and after by \emph{degree} we always mean
the Beauville-Bogomolov-Fujiki degree. The moduli space $\cV_{n,d}$ exists as a 
quasi-projective variety with quotient singularities (cf.\ \cite{Viehweg:ModuliSpaces}).

In \cite{Markman:Survey}, the moduli space
\[ \cV_{\bar{h}}:=\coprod \cV^t_{\bar{h}} \]
is referred to as the moduli space of polarised irreducible holomorphic symplectic manifolds of type $\bar{h}$.

Recall that the divisibility $\divv(h)$ of a non-zero and primitive element $h\in\Lambda$ is the positve generator of
the ideal $(h,\Lambda)\subset\ZZ$. Notice that this number always divides the determinant
of $\Lambda$. Fixing the degree and the divisibility of $h$ does not determine its orbit
$\bar{h}$ in general.
If $\cV_{n,d,\delta}\subset\cV_{n,d}$ is the sub-moduli space in which the polarisation
has divisibility $\delta$, then 
$$\cV_{n,d,\delta}\cong\coprod\cV_{\bar{h}},$$
where the disjoint union runs over all the $\Or(\Lambda)$-orbits of elements $h$ with degree
$d$ and divisibility $\delta$. 
In particular, if $\Upsilon_{n,d,\delta}$ is the set of connected components of
$\cV_{n,d,\delta}$, then
$$\Upsilon_{n,d,\delta}=\coprod\Upsilon_{\bar{h}}.$$
Notice that all these sets are finite.

The sub-moduli spaces $\cV_{n,d,\delta}$ are the main objects in 
Section~\ref{section:polarised}.

%%%%%%%%%%%%%%%%%%%%%%%%%%%%%%%%%%%%%%%%%
\section{Generalised Kummer varieties and a characterisation of parallel transport operators}\label{section:pto}

The main characters of this paper are the so-called generalised Kummer varieties. We recall
the two main constructions of such manifolds.
\begin{example}[\cite{Beauville:c1=0}]\label{exe:Kum}
Let $A$ be an abelian surface and $n\geq2$. The Hilbert scheme $A^{[n+1]}$ of 
$(n+1)$-points on $A$ is smooth and its Albanese map 
$a\colon A^{[n+1]}\to A$ is just the sum map (defined using the group structure of $A$).
If $p\in A$, then the fibre $a^{-1}(p)=:\Kum^n(A)$ is an irreducible holomorphic symplectic
manifold of dimension $2n$.
\end{example}
Let $A$ be an abelian surface. The even cohomology ring 
$H^{\operatorname{even}}(A,\ZZ)=H^0(A,\ZZ)\oplus H^2(A,\ZZ)\oplus H^4(A,\ZZ)$ has a 
natural quadratic form
$$ (\alpha_1,\alpha_2,\alpha_3)^2=\int_A(\alpha_2^2-2\alpha_1\alpha_3)$$
turning it into an even unimodular lattice called \emph{Mukai lattice}. 
It has signature $(4,4)$ and is isometric to the abstract
lattice $\widetilde{\Lambda}:=U^{\oplus4}$, where $U$ is the hyperbolic plane. 
\begin{example}[\cite{Yoshioka:ModuliSpacesAbelianSurfaces}]\label{exe:K}
Let again $A$ be an abelian surface and let $v\in H^{\operatorname{even}}(A,\ZZ)$ be
a primitive and effective Mukai vector such that $v^2\geq6$ 
(see \cite[Definition~0.1]{Yoshioka:ModuliSpacesAbelianSurfaces}). 
The moduli space $M_v(A)$ (with respect to a $v$-generic polarisation) is again smooth 
and the Albanese map is
$$\mathfrak{a}\colon M_v(A)\longrightarrow A\times A^\vee$$
defined by 
$$\mathfrak{a}(E)=\left(\det\left(\mathfrak{F}(E)\otimes\mathfrak{F}(E_0)^\vee\right), \det(E)\otimes\det(E_0)^\vee\right),$$
where $\mathfrak{F}$ is the Fourier-Mukai transform, $E_0\in M_v(A)$ is fixed and
$A^\vee=\Pic^0(A)$ is the dual abelian variety.
The fibre $\mathfrak{a}^{-1}(E_0)=:\Kt_v(A)$ is an irreducible holomorphic symplectic
manifold of dimension $v^2-2$ deformation equivalent to $\Kum^{\frac{v^2}{2}-1}(A)$.
\end{example}
Any irreducible holomorphic symplectic manifold deformation equivalent to one of the examples
above is called of \emph{Kummer type}. If $X$ is one such manifold, the
lattice structure on $H^2(X,\ZZ)$ is isometric to
\begin{equation}
\Lambda_n:=U^{\oplus3}\oplus\langle-2n-2\rangle,
\end{equation}
where $U$ is the hyperbolic plane.

The discriminant group $A_X$ of an irreducible holomorphic symplectic manifold is the quotient
$H^2(X,\ZZ)^*/H^2(X,\ZZ)$; any isometry $g\in\Or(H^2(X,\ZZ))$ naturally acts on $A_X$. 
Define 
\begin{equation}
W(X)=\{g\in\Or(H^2(X,\ZZ))\mid g \mbox{ acts as } \pm id \mbox{ on } A_X\}
\end{equation}
and consider the associated character
$
\chi\colon W(X)\longrightarrow\{\pm1\}.
$
Let $f\colon W(X)\longrightarrow\{\pm1\}$
be the map $f(g)=\chi(g)\det(g)$ and define
\begin{equation}
N(X)=\ker f.
\end{equation}
\begin{rmk}\label{rmk:W(X) as reflections}
If $u\in H^2(X,\ZZ)$ is such that $q(u)=\pm2$, define the reflections
\begin{equation}\label{appendix eqn:reflections}
\rho_u(v)=\left\{
\begin{array}{cc}
v+q(u,v)u & q(v)=-2 \\
-v+q(u,v)u & q(v)=2
\end{array}\right.
\end{equation}
Notice that $W(X)$ is the group generated by products of reflections $\rho_u$, 
where $(u,u)=\pm2$ (\cite[Lemma~4.2]{Markman:IntegralConstraints}). 
It follows that $N(X)$ is the group generated
by products $\rho_{u_1}\cdots\rho_{u_k}$, where $(u_j,u_j)=-2$ for an even number of 
indices, and $(u_j,u_j)=2$ for the remaining ones. 
\end{rmk}
\begin{prop}[\protect{\cite[Theorem~2.3]{Mongardi:Monodromy}}]\label{thm:monodromy groups}
Let $X$ be an irreducible holomorphic symplectic manifold of Kummer type. Then
$$\Mon^2(X)=N(X).$$
\end{prop}
When $X=\Kt_v(A)$ is a moduli space as in Example~\ref{exe:K}, 
we have a natural isometry $H^2(X,\ZZ)\cong v^\perp$ and so
a natural primitive embedding 
$i_v\colon H^2(X,\ZZ)\to H^{\operatorname{even}}(A,\ZZ)\cong\widetilde{\Lambda}$. 
\begin{rmk}\label{rmk:W(X)}
Notice that if $g\in W(X)$ then it extends to the lattice
$H^{\operatorname{even}}(A,\ZZ)$, i.e.\ there exists an isometry 
$\tilde{g}\in\Or(H^{\operatorname{even}}(A,\ZZ))$ such that 
$\tilde{g}|_{H^2(X,\ZZ)}=g$ (\cite[Proposition 1.5.1]{Nikulin:Lattice}). 
\end{rmk}

Let $\Or(\Lambda_n,\widetilde{\Lambda})$ be the set of primitive embeddings 
of $\Lambda_n$ inside $\widetilde{\Lambda}$. Both $\Or(\Lambda_n)$ and 
$\Or(\widetilde{\Lambda})$ act on $\Or(\Lambda_n,\widetilde{\Lambda})$ by, 
respectively, pre- and post-composition.
\vspace{0.3cm}

\begin{prop}[\protect{\cite[Theorem~4.9]{Wieneck:MonodromyInvariantsFibrations}}]\label{prop:i_X}
There exists a distinguished $\Mon^2(X)$-invariant $\Or(\widetilde{\Lambda})$-orbit
$$[i_X]\in\Or(\widetilde{\Lambda})\backslash\Or(H^2(X,\ZZ),\widetilde{\Lambda}).$$
\end{prop}
Let us recall how this orbit is constructed.
Deform $X$ to $\Kt_v(A)$ and pick
a parallel transport operator $P\colon H^2(X,\ZZ)\to H^2(\Kt_v(A),\ZZ)$.
As we said above, there exists a distinguished primitive embedding 
$i_v\colon H^2(K(v),\ZZ)\to H^{\operatorname{even}}(A,\ZZ)$ and hence a distinguished 
$\Or(\widetilde{\Lambda})$-orbit 
$[i_v]\in\Or(\widetilde{\Lambda})\backslash\Or(H^2(K(v),\ZZ),\widetilde{\Lambda})$.
Put then 
$$
[i_X]:=[i_v\circ P]\in\Or(\widetilde{\Lambda})\backslash\Or(H^2(X,\ZZ),\widetilde{\Lambda}).
$$
\begin{rmk}
Notice that $W(X)$ is identified with the 
stabiliser with respect to the $\Or^+(H^2(X,\ZZ))$-action of $[i_X]$ in $\Or(\widetilde{\Lambda})\backslash\Or(H^2(X,\ZZ),\widetilde{\Lambda})$ (\cite[Lemma~4.3]{Markman:IntegralConstraints}).
\end{rmk}

\begin{prop}\label{prop:characterisation pto}
Let $X_1$ and $X_2$ be two irreducible holomorphic symplectic manifolds of Kummer 
type, and let $g\colon H^2(X_1,\ZZ)\to H^2(X_2,\ZZ)$ be an orientation preserving isometry.
Then,
\begin{enumerate}
\item if $g$ is a parallel transport operator, then $[i_{X_1}]=[i_{X_2}]\circ g$;
\item if $[i_{X_1}]=[i_{X_2}]\circ g$, then either $g$ is a parallel transport operator
or $\tau_{X_2}\circ g$ is, where $\tau_{X_2}$ is any element in $W(X_2)\setminus N(X_2)$.
\end{enumerate}
\end{prop}
Notice that since $N(X)$ has index $2$ in $W(X)$, the choice of $\tau_{X_2}$ is essentially unique.
\begin{proof}
Assume first that $g$ is a parallel transport operator.
Let us deform both $X_1$ and $X_2$ to the same moduli space $\Kt_v(A)$ and pick two
parallel transport operators $P_i\colon H^2(X_i,\ZZ)\to H^2(\Kt_v(A),\ZZ)$. 
By assumption on $g$, we can choose $P_1=P_2\circ g$ and then
$$[i_{X_1}]=[i_v\circ P_1]=[i_v\circ P_2\circ g]=[i_v\circ P_2]\circ g=[i_{X_2}]\circ g$$ 
by the definition of $[i_X]$.

Vice versa, let us suppose that $[i_{X_1}]=[i_{X_2}]\circ g$. Since $X_1$ and $X_2$ are 
deformation equivalent, we can pick a parallel transport operator 
$f\colon H^2(X_2,\ZZ)\to H^2(X_1,\ZZ)$ and by the previous part of the proof we have
$[i_{X_2}]=[i_{X_1}]\circ f$.
Putting together these two equalities, we get the relation 
$[i_{X_1}]=[i_{X_1}]\circ(f\circ g)$, that is $f\circ g\in W(X_1)$.

If $f\circ g\in N(X_1)$, then we conclude as before.
If $f\circ g\notin N(X_1)$, then there exists $h\in W(X_1)\setminus N(X_1)$
such that $h\circ f\circ g\in N(X_1)$ is a monodromy operator. As before, the composition
$(f^{-1}\circ h\circ f)\circ g$ is a parallel transport operator and 
$f^{-1}\circ h\circ f=\tau_X$ is the required element in $W(X_1)$ but not in $N(X_1)$.
\end{proof}

Now let $(X_1,H_1)$ and $(X_2,H_2)$ be two polarised deformation equivalent irreducible
holomorphic symplectic manifolds. Let us put $h_i:=c_1(H_i)$ for convenience.
Using \cite[Proposition~7.4]{Markman:Survey}, we get the following corollary.
\vspace{0.3cm}

\begin{cor}\label{cor:characterisation ppto}
Suppose $X_1$ and $X_2$ are irreducible holomorphic symplectic manifolds of generalised 
Kummer type, and
let $g\colon H^2(X_1,\ZZ)\to H^2(X_2,\ZZ)$ be an orientation preserving isometry.
Then:
\begin{enumerate}
\item if $g$ is a polarised parallel transport operator, then $[i_{X_1}]=[i_{X_2}]\circ g$
and $g(h_1)=h_2$;
\item if $[i_{X_1}]=[i_{X_2}]\circ g$ and $g(h_1)=h_2$, then either $g$ is a polarised 
parallel transport operator or there exists an element $u\in H^2(X_2,\ZZ)$, with 
$(u,u)=-2$ and $(u,h_2)=0$, such that $\rho_u\circ g$ is a parallel transport operator.
\end{enumerate}
\end{cor}
\begin{proof}
By Remark~\ref{rmk:W(X) as reflections}, $\rho_u$ is an element in $W(X_2)$ but not
in $N(X_2)$ as soon as $(u,u)=-2$.
The only thing to prove is then the existence of elements $u$ such that $(u,u)=-2$ and 
$(u,h_2)=0$.
Since $(-2)$-elements exist in hyperbolic planes, this follows from the Eichler criterion 
(\cite[Proposition 3.3]{GritsenkoHulekSankaran:Abelianisation}).
\end{proof}

%%%%%%%%%%%%%%%%%%%%%%%%%%%%%%%%%%%%%%%%
\section{Moduli spaces of marked irreducible holomorphic symplectic manifolds}\label{section:marked}
Let $\Lambda_n=U^{\oplus3}\oplus\langle-2n-2\rangle$ be the abstract lattice of an 
irreducible holomorphic symplectic manifold of dimension $2n$ of Kummer type, and 
$\fM_{\Lambda_n}$ the moduli space of marked irreducible holomorphic symplectic manifolds of Kummer type.
We denote by $\Xi_n$ the set of connected components of $\fM_{\Lambda_n}$.

By Remark~\ref{rmk:orient pub}, any $X$ comes with a preferred orientation on $H^2(X,\ZZ)$ (induced by the choice of a K\"ahler class and a holomorphic symplectic class). If $\eta$ is a marking, then $\eta$ induces an orientation on the abstract lattice $\Lambda_n$. (As in Remark~\ref{rmk:orient pub}, one defines the orientation of a lattice $\Lambda$ of signature $(3,t)$ to be a generator of $H^2(\widetilde{C}_\Lambda,\ZZ)$, where $\widetilde{C}_\Lambda$ is the cone of positive vectors - see \cite[Section~4]{Markman:Survey}.) Let us call $\operatorname{orient}(X,\eta)$ this 
orientation.
As explained in \cite[Section~4]{Markman:Survey}, the induced orientation map
$$
\operatorname{orient}\colon\Xi_n\longrightarrow\operatorname{orient}(\Lambda_n)=\{\pm1\},
$$
defined by sending each connected component $\fM_{\Lambda_n}^t$ to 
$\operatorname{orient}(X,\eta)$ for any $(X,\eta)\in\fM_{\Lambda_n}^t$, is well defined 
and independent of the choice of $(X,\eta)$. Essentially this is due to the fact that if $(X_1,\eta_1)$ and $(X_2,\eta_2)$ are in the same connected component, then the composition $\eta_2^{-1}\circ\eta_1$ is a parallel transport operator, and hence it is orientation preserving.

On the other hand we can define the map
$$
\operatorname{orb}\colon\Xi_n\longrightarrow\Or(\widetilde{\Lambda})\backslash\Or(\Lambda_n,\widetilde{\Lambda})
$$
by sending $\fM_{\Lambda_n}^t$ to $[i_X]\circ\eta^{-1}$, where $(X,\eta)\in\fM_{\Lambda_n}^t$.
This is well defined, because if $(X',\eta')\in\fM_{\Lambda_n}^t$ is another marked pair, 
then the composition $\eta^{-1}\circ\eta'$ is a parallel transport operator and, by 
Proposition~\ref{prop:characterisation pto}, $[i_{X'}]\circ\eta'^{-1}=[i_X]\circ\eta^{-1}$.
\vspace{0.3cm}

\begin{prop}\label{prop:orb x orient}
The product map
$$
\operatorname{orb}\times\operatorname{orient}\colon\Xi_n\longrightarrow\Or(\widetilde{\Lambda})\backslash\Or(\Lambda_n,\widetilde{\Lambda})\times\{\pm1\}
$$
is 2:1 and surjective.
\end{prop}
\begin{proof}
It directly follows from Proposition~\ref{prop:characterisation pto}.
\end{proof}
\begin{cor}\label{cor:marked}
The number of connected components of the moduli space $\fM_{\Lambda_n}$ of marked irreducible
holomorphic symplectic manifolds of Kummer type is
$$|\Xi_n|=2^{\rho(n+1)+1},$$
where $\rho(k)$ is the number of distinct primes in the factorisation of $k$.
\end{cor}
\begin{proof}
By \cite[Lemma~4.3.(1)]{Markman:IntegralConstraints}, the cardinality of
$\Or(\widetilde{\Lambda})\backslash\Or(\Lambda_n,\widetilde{\Lambda})$
is equal to $2^{\rho(n+1)-1}$. 
Hence the claim follows directly from Proposition~\ref{prop:orb x orient}.
\end{proof}

%%%%%%%%%%%%%%%%%%%%%%%%%%%%%%%%%%%%%%%
\section{Moduli spaces of polarised irreducible holomorphic symplectic manifolds}\label{section:polarised}
Let $\Lambda_n$ be the abstract lattice of an irreducible holomorphic symplectic 
manifold of dimension $2n$ of Kummer 
type and let $h\in\Lambda_n$ be
a primitive element such that $(h,h)=2d>0$ (notice that $\Lambda_n$ is even). Denote by 
$\bar{h}$ the $\Or(\Lambda_n)$-orbit of $h$.
Recall from Section~\ref{section:pre} that 
$\mathcal{V}^t_{\bar{h}}=\fM_{\bar{h}}/\Or(\Lambda_n)$ 
is (isomorphic to) one connected component of the moduli space of polarised irreducible holomorphic symplectic 
manifolds $(X,H)$ of dimension $2n$ and such that $q_X(c_1(H))=2d$.
If $\delta$ is a positive divisor of $2n+2$, then $\cV_{n,d,\delta}\cong\coprod\cV^t_{\bar{h}}$ 
is (isomorphic to) the moduli space of polarised manifolds with polarisation of degree $2d$ and divisibility $\delta$. Here the sum runs over all the $\Or(\Lambda_n)$-orbits of vectors of degree $2d$ and divisibility $\delta$, and over all the connected components of the moduli space $\cV_{\bar{h}}$ (see Section~\ref{section:pre}).
We want to compute the cardinality of the set $\Upsilon_{n,d,\delta}$ of conneceted components
of $\cV_{n,d,\delta}$.

Let $(X,H)\in\fM_{\bar{h}}/\Or(\Lambda_n)$ be a polarised pair and pick a representative 
$i\in[i_X]$.
The orthogonal complement $i(H^2(X,\ZZ))^\perp\subset\widetilde{\Lambda}$ is a positive rank 
$1$ sublattice. 
If $T_{(X,H)}$ is the saturation of the lattice generated by $i(H^2(X,\ZZ))^\perp$ and $i(c_1(H))$, then $T_{(X,H)}$ is a positive and primitive rank $2$ sublattice of $\widetilde{\Lambda}$.
%Therefore the lattice $T_{(X,H)}$, primitively generated by $i(H^2(X,\ZZ))^\perp$ and $i(c_1(H))$, is a positive rank $2$ sublattice of $\widetilde{\Lambda}$.

Notice that if $i'$ is another representative of $[i_X]$, then there exists an isometry 
$\tilde{g}\in\Or(\widetilde{\Lambda})$ which restricts to an isometry $g\in\Or(T_{(X,h)})$.
Moreover, by construction $g(i(c_1(H)))=i'(c_1(H))$.

This suggests the definition of the following set
$$\Sigma_{n}=\left\{(T,h)\mid T \mbox{ positive rank $2$ lattice and } h\in T \mbox{ primitive s.t.\ } h^\perp=\langle2n+2\rangle\right\}/\sim,$$
where $(T,h)\sim(T',h')$ if there exists an isometry $g\colon T\to T'$ such that $g(h)=h'$.
We denote by $[T,h]$ the equivalence classes. 
\begin{rmk}\label{rmk:Nik}
By \cite[Theorem 1.1.2]{Nikulin:Lattice}, $T$ can be primitively 
embedded in $\widetilde{\Lambda}$ in a unique way (up to an isometry of $\widetilde{\Lambda}$).
\end{rmk}

Now let $I(X)$ be the set of positive and primitive classes in $H^2(X,\ZZ)$. 
There is a well-defined map
$$ f_X\colon I(X)\longrightarrow\Sigma_n $$
defined by sending $h\in I(X)$ to $[T_{(X,h)},i(h)]$ for any $i\in[i_X]$.
In the following, we drop the dependence on $[i_X]$ from the notation and we simply write $[T{(X,h)},h]$.
\vspace{0.3cm}

\begin{prop}\label{prop:f for ppto}
Given two polarised pairs $(X_1,H_1)$ and $(X_2,H_2)$ of manifolds of Kummer 
type, a polarised parallel transport operator
between them exists if and only if $f_{X_1}(c_1(H_1))=f_{X_2}(c_1(H_2))$.
\end{prop}
\vspace{0.3cm}
The following proof is a translation to our case of the proof of 
\cite[Proposition~1.6]{Apostolov:ConnectedComponents}.
For sake of notation, we put $h_i=c_1(H_1)$. 
\begin{proof}
Suppose that $P\colon H^2(X_1,\ZZ)\to H^2(X_2,\ZZ)$ is a polarised parallel transport operator.
By \cite[Proposition~7.4]{Markman:Survey} and Corollary~\ref{cor:characterisation ppto}, we immediately get
an isometry $\tilde{g}\in\Or(\widetilde{\Lambda})$ which restricts to an isometry 
$g\colon T_{(X_1,h_1)}\to T_{(X_2,h_2)}$ such that $g(h_1)=h_2$.

Vice versa, suppose that such an isometry $g$ exists. In particular $T_{(X_1,h_1)}$ has two primitive
embeddings inside $\widetilde{\Lambda}$, the second one given by composing the natural embedding
$T_{(X_2,h_2)}\subset\widetilde{\Lambda}$ with $g$. 
By Remark~\ref{rmk:Nik}, there exists a unique (up to isometry) such
primitive embedding and hence there exists an isometry $\tilde{g}\in\Or(\widetilde{\Lambda})$ such that
the diagram
\begin{equation*}
\xymatrix{
\widetilde{\Lambda}\ar@{->}[r]^{\tilde{g}} & \widetilde{\Lambda} \\
T_{(X_1,h_1)}\ar@{^{(}->}[u]\ar@{->}[r]^{g} & T_{(X_2,h_2)}\ar@{^{(}->}[u]
}
\end{equation*}
commutes. 

Since $\tilde{g}(i_1(H^2(X_1,\ZZ)))=i_2(H^2(X_2,\ZZ))$, it follows that $\tilde{g}$
restricts to an isometry $P$ from $H^2(X_1,h_1)$ to $H^2(X_2,h_2)$. Here $i_1\in[i_{X_1}]$
and $i_2\in[i_{X_2}]$, and we have a commutative diagram
\begin{equation*}
\xymatrix{
\widetilde{\Lambda}\ar@{->}[r]^{\tilde{g}} & \widetilde{\Lambda} \\
H^2(X_1,\ZZ)\ar@{^{(}->}[u]^{i_1}\ar@{->}[r]^{P} & H^2(X_2,\ZZ)\ar@{^{(}->}[u]^{i_2}.
}
\end{equation*}
In particular $[i_{X_2}]\circ P=[i_{X_1}]$ and $P(h_1)=h_2$.
By Corollary~\ref{cor:characterisation ppto}, we get a polarised parallel transport operator from $P$
as long as $P$ is orientation preserving.

Let us then suppose that $P$ is not orientation preserving. Let us pick an element $u\in H^2(X_2,\ZZ)$
such that $(u,u)=2$ and $(u,h_2)=0$, and let us consider the reflection $\rho_u$. Since $\rho_u$
is orientation preserving by definition and $\rho_u(h_2)=-h_2$, then $P'=-\rho_u\circ P$ is an orientation 
preserving isometry such that $[i_{X_2}]\circ P'=[i_{X_2}]$ and $P'(h_2)=h_2$, and we can
apply Corollary~\ref{cor:characterisation ppto} to produce a polarised parallel transport operator.
\end{proof}

\begin{rmk}
$f_X$ is a \emph{faithful monodromy invariant}, as defined in \cite[Section~5.3]{Markman:PrimeExceptional}.
\end{rmk}
For the next result, we define $\cV_n:=\coprod\cV_{\bar{h}}$, where the disjoint union runs
over all the $\Or(\Lambda)$-orbits (in particular, we are not even fixing the degree of the 
polarisation). The set of the connected components of $\cV_n$ is denoted by $\Upsilon_n=
\coprod\Upsilon_{\bar{h}}$. Furthermore, define $\Sigma_{n,d,\delta}\subset\Sigma_{n}$ as 
the subset consisting of pairs $[T,h]$ such that $(h,h)=2d$ and $\divv(h)=\delta$; 
notice that $\Sigma_n=\coprod\Sigma_{n,d,\delta}$.

The existence of a polarised parallel transport operator between $(X_1,H_1)$ and $(X_2,H_2)$
is equivalent to saying that $(X_1,H_1)$ and $(X_2,H_2)$ belong to the same connected component.
\vspace{0.3cm}

\begin{prop}\label{prop:f}
There exists a well-defined injective map
\begin{equation}\label{eqn:f polarised}
f\colon \Upsilon_n\longrightarrow\Sigma_n
\end{equation}
defined by sending a connected component $\cV^t_{\bar{h}}$ to $f_{X}(c_1(H))$, for any 
$(X,H)\in\cV^t_{\bar{h}}$.

Moreover, $f(\Upsilon_{n,d,\delta})=\Sigma_{n,d,\delta}$.
\end{prop}

\begin{proof}
The proof is the same as the proof of 
\cite[Theorem~1.7, Proposition~2.3]{Apostolov:ConnectedComponents}, up to use 
Proposition~\ref{prop:f for ppto} above istead of 
\cite[Proposition~1.6]{Apostolov:ConnectedComponents}.
\end{proof}

In the rest of this section we want to compute the cardinality of $\Sigma_{n,d,\delta}$. 
The first remark is that a pair $[T,h]\in\Sigma_n$ is completely determined by the primitive embedding
$j\colon\langle h\rangle\to T$ such that $j(h)^\perp=\langle2n+2\rangle$.
Therefore we want to count the number of such primitive embeddings.
Let us note that, by \cite[Theorem~1.1.2]{Nikulin:Lattice}, without loss of generality we can
think of both $\langle h\rangle$ and $T$ as sublattices of $\widetilde{\Lambda}$.

The main result of this section is the following. 
First of all, we make the following definitions (which will be useful during the proof of the 
theorem below):
\begin{equation}\label{eqn:tante}
\begin{array}{lll}
d_1=\frac{2d}{\gcd(2d,2n+2)}, & n_1=\frac{2n+2}{\gcd(2d,2n+2)}, & g=\frac{\gcd(2d,2n+2)}{\delta} \\
 & & \\
w=\gcd(g,\delta), & g_1=\frac{g}{w}, & \delta_1=\frac{\delta}{w}.
\end{array}
\end{equation}
Furthermore, the following notation is used in the statement of the main theorem:
for an integer $l$ we write $\phi(l)$ for the Euler function and $\rho(l)$ for the number of distinct primes in the factorisation of $l$; for $w$ and $\delta_1$ as defined above, we write $w=w_+(\delta_1)w_-(\delta_1)$, where $w_+(\delta_1)$ is the product of all powers of the primes that appear in the factorisation of $w$ and that divide $\gcd(w,\delta_1)$ (that is, $w_-(\delta_1)$ is the part coprime to $\delta_1$). More precisely, if $p^k$ is a factor of $w$ and $p$ divides $\gcd(w,\delta_1)$, then $p^k$ is a factor of $w_+(\delta_1)$.
\vspace{0.3cm}

\begin{thm}\label{thm:connected components polarised}
With the notations as above, we have:
\begin{enumerate}

\item $|\Upsilon_{n,d,\delta}|=w_+(\delta_1)\phi(w_-(\delta_1))2^{\rho(\delta_1)-1}$ if $\delta>2$ and one of the following holds:
 \begin{enumerate}
 \item $g_1$ is even, $\gcd(d_1,\delta_1)=1=\gcd(n_1,\delta_1)$ and $-d_1/n_1$ is a quadratic residue mod~$\delta_1$;
 \item $g_1,\delta_1$ and $d_1$ are odd, $\gcd(d_1,\delta_1)=1=\gcd(n_1,2\delta_1)$ and $-d_1/n_1$ is a quadratic residue mod~$2\delta_1$;
 \item $g_1,\delta_1$ and $w$ are odd, $d_1$ is even, $\gcd(d_1,\delta_1)=1=\gcd(n_1,2\delta_1)$ and $-d_1/4n_1$ is a quadratic residue mod~$\delta_1$.
 \end{enumerate}

\item $|\Upsilon_{n,d,\delta}|=w_+(\delta_1)\phi(w_-(\delta_1))2^{\rho(\delta_1/2)-1}$ if $\delta>2$, $g_1$ is odd, $\delta_1$ is even, $\gcd(d_1,\delta_1)=1=\gcd(n_1,2\delta_1)$ and $-d_1/n_1$ is a quadratic residue mod~$2\delta_1$.

\item $|\Upsilon_{n,d,\delta}|=1$ if $\delta\leq2$ and one of the following holds:
\begin{enumerate}
 \item $g_1$ is even, $\gcd(d_1,\delta_1)=1=\gcd(n_1,\delta_1)$ and $-d_1/n_1$ is a quadratic residue mod~$\delta_1$;
 \item $g_1,\delta_1$ and $d_1$ are odd, $\gcd(d_1,\delta_1)=1=\gcd(n_1,2\delta_1)$ and $-d_1/n_1$ is a quadratic residue mod~$2\delta_1$;
 \item $g_1,\delta_1$ and $w$ are odd, $d_1$ is even, $\gcd(d_1,\delta_1)=1=\gcd(n_1,2\delta_1)$ and $-d_1/4n_1$ is a quadratic residue mod~$\delta_1$;
 \item $g_1$ is odd, $\delta_1$ is even, $\gcd(d_1,\delta_1)=1=\gcd(n_1,2\delta_1)$ and $-d_1/n_1$ is a quadratic residue mod~$2\delta_1$.
 \end{enumerate}

\item $|\Upsilon_{n,d,\delta}|=0$ otherwise.
\end{enumerate}
\end{thm}

\begin{proof}
Using the bijection (\ref{eqn:f polarised}) and the discussion above, $|\Upsilon_{n,d,\delta}|=|\Sigma_{n,d,\delta}|$ and the latter is the number of primitive embeddings
$j\colon\langle2d\rangle\to T$ such that $j(\langle2d\rangle)^\perp=\langle2n+2\rangle$.

By \cite[Proposition~1.5.1]{Nikulin:Lattice}, an embedding $j\colon\langle2d\rangle\to T$ is determined
by the pair $(H,\gamma)$, where $H\subset A_{2d}$ is a subgroup, 
$\gamma\colon H\to A_{2n+2}$ is an injective homomorphism and the pushout 
$\Gamma_\gamma=H\subset A_{2d}\oplus A_{2n+2}$ is isotropic.
Since we have also fixed $\divv(h)=\delta$, it follows that $H$ must be of order $\delta$ 
(see \cite[Proposition~2.2]{Apostolov:ConnectedComponents}).
\begin{rmk}\label{rmk:2d}
Recall  that two pairs $(H,\gamma)$ and $(H',\gamma')$ determine the same primitive embedding $j$ 
if $H=H'$ and there exist an isometry $\varphi\in\Or(\langle2d\rangle)\cong\ZZ/2\ZZ$ and an isometry
$\psi\in\Or(\langle2n+2\rangle)\cong\ZZ/2\ZZ$ such that
$\gamma\circ\overline{\varphi}=\overline{\psi}\circ\gamma'$ (\cite[Section~5]{Nikulin:Lattice}).
\end{rmk}
Identifying $A_{2d}$ with $\ZZ/2d\ZZ$ and picking generators $h$ of $\langle2d\rangle$ and $v$ 
of $\langle2n+2\rangle$, we can write $H=\langle h/\delta\rangle$. 
Then $\gamma$ is uniquely determined by the image $\gamma(h/\delta)=cv/\delta$, where $c$ 
is coprime with $\delta$. The isotropy condition is
\begin{equation}\label{eqn:1}
\frac{2d}{\delta^2}+\frac{c^2(2n+2)}{\delta^2}\equiv 0\pmod{2}.
\end{equation}
Substituting (\ref{eqn:tante}) in equation (\ref{eqn:1}), we eventually get
\begin{equation}\label{eqn:c}
\delta_1\left(\frac{2d}{\delta^2}+\frac{c^2(2n+2)}{\delta^2}\right)=g_1(d_1+c^2n_1)\equiv0\pmod{2\delta_1}.
\end{equation}

The problem is now reduced to determine all the solutions $c$ of equation (\ref{eqn:c}) 
such that $\gcd(c,\delta)=1$. This problem has already been solved by Gritsenko, Hulek and Sankaran in the proof of
\cite[Proposition~3.6]{GritsenkoHulekSankaran:Moduli}. Since we are interested in isometric embeddings,
we have to understand which of these solutions are invariant under the isometries in Remark~\ref{rmk:2d}.
Both $\Or(\langle2d\rangle)$ and $\Or(\langle2n+2\rangle)$ act on $H$ by changing the sign
of the first, respectively the second, coordinate. Moreover, notice that $H$ has a central symmetry, i.e.\ 
$(x,y)\in H$ if and only if $(-x,-y)\in H$. We can then distinguish two cases:
\begin{itemize}
\item $\delta\leq2$: then any subgroup $H$ is fixed by this action and the number of solutions $c$ corresponds to the 
number of primitive embeddings;
\item $\delta>2$: then there are no fixed subgroups $H$ and we must divide the number of solutions $c$ by $2$.
\end{itemize}
This concludes the proof.
\end{proof}

\begin{rmk}\label{rmk:w 1}
When $w=1$, the values of $d$ and $\delta$ determine the orbit of $h$, i.e.\ $\cV_{n,d,\delta}=\cV_{\bar{h}}$
(cf.\ \cite[Corollary 3.7]{GritsenkoHulekSankaran:Moduli}). 
\end{rmk}
We conclude this section by giving a few examples.
\begin{example}\label{example:delta 1}
If $\delta=1$, then the orbit of $h$ is determined and moreover the corresponding moduli space is connected.
\end{example}
\begin{example}
Let $p$ and $q$ be two (different) odd primes and put $\delta=d=pq$ and $n+1=mpq$, where $\gcd(m,pq)=1$ and $-m$ is a quadratic residue mod~$pq$. 
Then the moduli space $\mathcal{V}_{n,d,\delta}$ has two connected components. 
\end{example}
\begin{example}\label{example:square free}
If $\gcd(2d,2n+2)$ is square free, then $w=1$ 
(cf.\ \cite[Remark~3.13]{GritsenkoHulekSankaran:Moduli}). 
This is the case, for example, when $2n+2$ is square free.
\end{example}
\begin{example}
Let $X$ be a fourfold of Kummer type. Then $2n+2=6$ is square free, so that $w=1$ and the values of the degree and the divisibility determine the $\Or(\Lambda_2)$-orbit of the polarisation (see Example~\ref{example:square free} and Remark~\ref{rmk:w 1}). We claim that the moduli space $\cV_{2,d,\delta}=\cV_{\bar{h}}$ is always connected.

In fact suppose that $h=c_1(H)$ is the class of a polarisation on $X$ with $q_X(h)=2d$ and $\divv(h)=\delta$.

\begin{enumerate}
\item If $\delta=1$, then this is Example~\ref{example:delta 1}.
\item If $\delta=2$, then $2d=8k-6$ for some integer $k$. We have two cases: either $k$ is a multiple of $3$, or $k$ is coprime to $3$. In the first case $d_1=4k'-1$, $n_1=1$, $g_1=3$ and $-d_1/n_1=1-4k'\equiv1$ mod~$4$. In the second case $d_1=4k-3$, $n_1=3$, $g_1=1$ and $-d_1/n_1=9-12k\equiv1$ mod~$4$. In both cases we have just one connected component by Theorem~\ref{thm:connected components polarised}.(3.(d)).
\item If $\delta=3$, then $2d=18k-6t^2$ for some integer $t$ coprime to $3$. Then $d_1=3k-t^2$, $n_1=1$, $g_1=2$ and $-d_1/n_1\equiv t^2$ mod~$3$ is a quadratic residue. Then the claim follows from Theorem~\ref{thm:connected components polarised}.(1.(a)).
\item If $\delta=6$, then $2d=72k-6t^2$ for some integer $t$ coprime to $6$. Then $d_1=12k-t^2$, $n_1=1$, $g_1=1$ and $-d_1/n_1\equiv t^2$ mod~$12$ is a quadratic residue. The claim follows from Theorem~\ref{thm:connected components polarised}.(2).
\end{enumerate}

\end{example}

\begin{example}
Let $X$ be a manifold of Kummer type of dimension $6$, and suppose that $h$ is the class of a polarisation such that $\divv(h)=4$. The degree $q_X(h)=2d$ must be of the form $2d=32k-8t^2$, with $t$ coprime to $4$ and one can check that $g=w=2$. Therefore we are in the situation of Theorem~\ref{thm:connected components polarised}.(2) and the number of conncted components of $\cV_{3,d,4}$ is $2$.
\end{example}

To conclude, let us remark that the number of connected components can get arbitrarily
large as the dimension and the degree of the polarisation increase.
\vspace{0.5cm}

%%%%%%%%%%%%%%%% BIBLIOGRAPHY %%%%%%%%%%%%%%%

\bibliographystyle{alpha}
% here we change the meaning of \VAN to use the prefix for the bibliography
\DeclareRobustCommand{\VAN}[3]{#3}
\bibliography{Bibliography}
\addcontentsline{toc}{chapter}{Bibliography}
%%%%%%%%%%%%%%%%%%%%%%%%%%%%%%%%%%%%%%%
%\affiliationone{% in this example, two authors share an institution
  % Claudio Onorati\\
  % Dipartimento di matematica \\
  % Universit\'a di Bologna \\
   %   Piazza di Porta San Donato, 5\\
    %  40126 Bologna, BO, Italia
   % \email{claudio.onorati@unibo.it}}

\end{document}